\documentclass[reqno]{amsart}

\usepackage{amsmath}
\usepackage{amsfonts}
\usepackage{amssymb}
\usepackage{verbatim}

\newtheorem{theorem}{Theorem}[section]

\theoremstyle{definition}
\newtheorem{definition}[theorem]{Definition}

\theoremstyle{remark}
\newtheorem{remark}[theorem]{Remark}

\numberwithin{equation}{section}
\usepackage[pdftex]{hyperref}
\def\e{{\mathrm{e}}}

\begin{document}

\title[Noncommutative tori]{Heisenberg modules over real multiplication noncommutative tori and related algebraic structures}

\author{Jorge Plazas}
\address{ I.H.E.S.
Le Bois-Marie,
35 route de Chartres
91440 Bures-sur-Yvette,
France}
\email{plazas@ihes.fr}

\subjclass[2000]{Primary 58B34, 46L87; Secondary 16S38}
\keywords{noncommutative tori, real multiplication, Heisenberg groups}

\begin{abstract}

We review some aspects of the theory of noncommutative two-tori with real multiplication focusing on 
the role played by Heisenberg groups in the definition of algebraic structures associated to these noncommutative spaces. 

\end{abstract}

\maketitle
\setcounter{section}{-1}
\section{Introduction}
\label{intro}

Noncommutative tori have played a central role in noncommutative geometry since the early stages of the theory. They arise naturally in various context and have provided a good testing ground for many of the techniques from which noncommutative geometry has developed (\cite{Connes1, rieffel2}). Noncommutative tori are defined in terms of their algebras of functions. The study of projective 
modules over these algebras and the corresponding theory of Morita equivalences leads to the 
existence of a class of noncommutative tori related to real quadratic extensions of $\mathbb{Q}$. These real multiplication noncommutative tori are conjectured to provide the correct geometric 
setting under which to attack the explicit class field theory problem for real quadratic fields (\cite{Manin1}). The right understanding of the algebraic structures underlying  these spaces is important for these applications.

The study of connections on vector bundles over noncommutative tori gives rise to a rich theory  which has been recast recently in the context of complex algebraic geometry (\cite{Connes1, ConnesRieffel, Schwarz1,DiS, PoS}). The study of categories of holomorphic bundles has thrown light about some algebraic structures related to real multiplication noncommutative tori (\cite{Po, Masha, yo}). Some of these results arise in a natural way from the interplay between Heisenberg groups and noncommutative tori. 

I want to thank the organizers of the 2007 summer school ``Geometric and topological methods for quantum field theory" and I.H.E.S.  for their support and hospitality. This work was supported in part by  ANR-Galois grant NT05-2 44266.

\section{Noncommutative tori and their morphisms}

In many situations arising in various geometric settings  it is possible to characterize spaces and  some of their structural properties in terms of appropriate rings of functions. One instance of this duality is provided by Gelfand's theorem which identifies the category of locally compact Hausdorff topological spaces with the category commutative $C^{*}$-algebras. This correspondence assigns to a space $X$ the commutative  algebra $C_{0}(X)$ consisting of complex valued continuous functions on $X$ vanishing at infinity.  Topological invariants of the space $X$ can be obtained by the corresponding invariants of  $C_{0}(X)$ defined in the context of  $C^{*}$-algebras. If $X$ is a smooth manifold the smooth structure on $X$ singles out the $*$-subalgebra $C^{\infty}_{0}(X)$ consisting of smooth elements of $C_{0}(X)$. Considering  the space $X$ in the framework of differential topology leads to structures defined in terms of the algebra $C^{\infty}_{0}(X)$.


Various geometric notions which can be defined in terms of rings of functions on a space do not depend on the fact that the rings under consideration are commutative and can therefore be extended in order to consider noncommutative rings and algebras. In noncommutative geometry spaces are defined in terms of their rings of functions which are noncommutative analogs of commutative rings of functions. This passage is far from being just a translation of classical ideas to a noncommutative setting. Many extremely rich new  phenomena arises in this context (see \cite{Connes3}). The noncommutative setting also enriches the classical picture since noncommutative rings may  arise in a natural way from classical geometric considerations.

We will be considering noncommutative analogs of the two-torus  {$\mathbb{T}^2=S^{1}\times S^{1}$}. \linebreak The reader may consult \cite{Effros, Connes1, Connes3, Davidson, riff3, rieffel2, Manin1} for the proofs of the results on noncommutative tori mentioned in this section.

\subsection{The $C^{*}$-algebra $A_{\theta}$}
Under the Gelfand correspondence compact spaces correspond to commutative unital $C^{*}$-algebras and  $\mathbb{T}^2$ is dual to $C(\mathbb{T}^2)$. At a topological level a noncommutative two-torus is defined in terms of a unital noncommutative $C^{*}$-algebra which plays the role of its algebra of continuous functions. 

Given $\theta \in \mathbb{R}$ let $A_{\theta}$ be the universal $C^{*}$-algebra generated by two unitaries $U$ and $V$ subject to the relation:
\begin{eqnarray*}
UV=e^{2\pi \imath \theta}VU.
\end{eqnarray*}
Then
\begin{enumerate}
\item If $\theta \in \mathbb{Z}$ the algebra $A_{\theta}$ is isomorphic to $C (\mathbb{T}^2)$.
\item If $\theta \in \mathbb{Q}$ the algebra $A_{\theta}$ is isomorphic to the algebra of global sections of the endomorphism bundle of a complex vector bundle over  $\mathbb{T}^2$. 
\item If $\theta \in \mathbb{R} \setminus \mathbb{Q}$  the algebra $A_{\theta}$ is a simple $C^{*}$-algebra.
\end{enumerate}
For irrational values of $\theta$  we will refer  to $A_{\theta}$ as {\it the algebra of continuous  functions on the noncommutative torus $\mathbb{T}^2_\theta$}.  Thus, as a topological space the noncommutative torus $\mathbb{T}^2_\theta$ is defined as the dual object of $C (\mathbb{T}^2_{\theta}) := A_{\theta}$. 

There is a natural continuous action of the compact group $\mathbb{T}^2$ on the algebra $A_{\theta}$. This action can be given in terms of the generators $U$ and $V$ by:
\begin{eqnarray*}
\alpha_{\varphi}(U) &=& e^{2 \pi \imath \varphi_1}U  \\
\alpha_{\varphi}(V) &=& e^{2 \pi \imath \varphi_2}V
\end{eqnarray*}
where $ \varphi = (\varphi_1,\varphi_2) \in \mathbb{T}^2$.

One of the main structural properties of the algebra $A_{\theta}$  is the existence of a canonical trace whose values at each element is given by the average over $\mathbb{T}^2$ of the above action.
\begin{theorem}
Let $\theta$ be an irrational number. Then there exist a unique normalized 
trace:
\begin{eqnarray*}
\label{trace1}
\chi: A_{\theta} \rightarrow \mathbb{C}
\end{eqnarray*}
invariant under the action of $\mathbb{T}^2$.
\end{theorem}
For the remaining part of the article $\theta$ will denote an irrational number. Also, for any complex number $z\in \mathbb{C}$ we will use the notation:
\begin{eqnarray*}
 \e(z) = \exp (2 \pi \imath z), && \bar{ \e}(z) = \exp (-2 \pi \imath z)\, .
\end{eqnarray*}

\subsection{Smooth elements}

The above action of $\mathbb{T}^2$ on $A_{\theta}$ induces a smooth structure on the noncommutative torus $\mathbb{T}^2_{\theta}$. An element $a\in A_{\theta}$ is called smooth if the map 
\begin{eqnarray*}
\mathbb{T}^2  & \longrightarrow & A_{\theta}  \\
\varphi & \longmapsto& \alpha_{\varphi}(a)
\end{eqnarray*}
is smooth. The set of smooth elements of $A_{\theta}$ is a dense $*$-subalgebra which we denote by
$\mathcal{A}_{\theta}$. Elements in this subalgebra should be thought as smooth functions on the noncommutative torus  $\mathbb{T}^2_{\theta}$ thus we take $C^{\infty}(\mathbb{T}^2_{\theta}) := \mathcal{A}_{\theta}$. The algebra $\mathcal{A}_{\theta}$ can be characterized in the following way:
\begin{eqnarray*}
\mathcal{A}_\theta
       &=& \bigg\{ \sum_{n,m \in \mathbb{Z}} a_{n,m}U^nV^m \, \in \, A_{\theta} \; | \; \{a_{n,m} \} \in \mathcal{S}(\mathbb{Z}^2)\bigg\} 
\end{eqnarray*}
where $\mathcal{S}(\mathbb{Z}^2)$ denotes the space of sequences of rapid decay in $\mathbb{Z}^2$. 

In the algebra $\mathcal{A}_\theta$ the trace $\chi$  is given  by
\begin{eqnarray*}
\chi (\sum a_{n,m}U^nV^m ) &=& a_{0,0}
\end{eqnarray*}

The Lie algebra $L= \mathbb{R}^2$ of $\mathbb{T}^2$ acts on $\mathcal{A}_\theta$ by derivations. A basis for this action is given by the derivations: 
\begin{eqnarray*} 
\label{derivation1} 
\delta_{1}(U)= 2\pi \imath U; && \delta_{1}(V)= 0 \\
\label{derivation2} 
\delta_{2}(U)= 0; && \delta_{2}(V)= 2\pi \imath V.
\end{eqnarray*}

A complex parameter $\tau \in  \mathbb{C} \setminus \mathbb{R}$ induces a complex structure on 
$L = \mathbb{R}^2$ given by the isomorphism:
\begin{eqnarray*}
\mathbb{R}^2  & \longrightarrow &  \mathbb{C} \\
x = (x_{1}, x_{2}) &\longmapsto& \tilde{x} = \tau x_{1} + x_{2}. 
\end{eqnarray*}
The corresponding complex structure on $\mathcal{A}_\theta$ is given by the derivation:
\begin{eqnarray*}
\delta_{\tau} &=& \tau \delta_{1} + \delta_{2}  \, .
\end{eqnarray*}

\subsection{Vector bundles and $K$-theory}
If  $X$ is a smooth compact manifold the space of smooth sections of a vector bundle over $X$ is a finite-type projective module over $C^{\infty}(X)$ and any such a module arises in this way. In our setting finite-type projective right $\mathcal{A}_\theta$-modules will play the role of vector bundles over the noncommutative torus $\mathbb{T}^2_{\theta}$. 

We denote as above by $L = \mathbb{R}^2$ the Lie algebra of $\mathbb{T}^2$ acting as an algebra of derivations on $\mathcal{A}_\theta$. If $P$ is a finite-type projective right $\mathcal{A}_\theta$-module a connection on $P$ is given by an operator:
\begin{eqnarray*}
\nabla: P \rightarrow P\otimes L^{*}
\end{eqnarray*}
such that 
\begin{eqnarray*}
\nabla_{X}(\xi a) &=&   \nabla_{X}(\xi) \, a +    \xi \delta_{X}  a
\end{eqnarray*}
for all $X\in L, \xi\in P$ and  $a\in\mathcal{A}_\theta$.

The connection $\nabla$ is determined the operators 
\begin{eqnarray*}
\nabla_{i}: P \rightarrow P, \qquad i=1,2
\end{eqnarray*}
giving its values on the basis elements $\delta_{1}, \delta_{2}$ of $L$.

The $K_{0}$ group of $\mathcal{A}_\theta$ is by definition the enveloping group of the abelian semigroup given by isomorphism classes of finite-type projective right $\mathcal{A}_\theta$-modules together with direct sum. The trace $\chi$ extends to an injective morphism 
\begin{eqnarray*}
\mathrm{rk}: K_{0}(\mathcal{A}_\theta) \rightarrow \mathbb{R}
\end{eqnarray*}
whose image is 
\begin{eqnarray*}
\Gamma_\theta &=&\mathbb{Z} \oplus \theta \mathbb{Z} \, .
\end{eqnarray*}

\subsection{Morphisms of noncommutative tori}

Since noncommutative tori are defined in terms of their function algebras one should expect a morphism $\mathbb{T}^2_{\theta} \rightarrow \mathbb{T}^2_{\theta'}$ between two noncommutative tori
$\mathbb{T}^2_{\theta}$ and $\mathbb{T}^2_{\theta'}$  to be given by a morphism 
$\mathcal{A}_{\theta'} \rightarrow \mathcal{A}_\theta$ of the corresponding algebras of functions. It turns out that algebra morphisms are in general insufficient to describe the type of situations arising in noncommutative geometry. The right notion of morphisms in our setting is given by Morita equivalences. 
A Morita equivalence between $\mathcal{A}_{\theta'}$ and $\mathcal{A}_{\theta}$ is given by the isomorphism class of a $\mathcal{A}_{\theta'}$-$\mathcal{A}_{\theta}$-bimodule 
$E$ which is projective and of finite-type both as a left \linebreak $\mathcal{A}_{\theta'}$-module and as a right $\mathcal{A}_{\theta}$-module. If such bimodule exists we say that $\mathcal{A}_{\theta'}$ and $\mathcal{A}_{\theta}$ are Morita equivalent. We can consider a Morita equivalence between $\mathcal{A}_{\theta'}$ and $\mathcal{A}_{\theta}$ as a morphism between 
$\mathcal{A}_{\theta'}$ and $\mathcal{A}_{\theta}$ inducing a morphism between $\mathbb{T}^2_{\theta}$ and $\mathbb{T}^2_{\theta'}$. Composition of morphisms is provided by tensor product of modules.

Let $SL_2 (\mathbb{Z})$ act on $\mathbb{R}\setminus \mathbb{Q}$ by fractional linear transformations, i.e.  given 
\begin{eqnarray*}
g= \left(
\begin{array}{cc}
  a & b \\
  c & d \\
\end{array}
\right)\in SL_2 (\mathbb{Z}), \qquad \theta \in \mathbb{R} \setminus \mathbb{Q}
\end{eqnarray*}
we take
\begin{eqnarray*}
g \theta &=& \frac{a\theta +b}{c\theta + d}\; .
\end{eqnarray*}
Morita equivalences between noncommutative tori are characterized by the following result:  
\begin{theorem}\emph{(Rieffel \cite{rieffel1})}
\label{Morita}
 Let $\theta', \, \theta \, \in \mathbb{R} \setminus \mathbb{Q}$. Then the algebras $\mathcal{A}_{\theta'}$ and $\mathcal{A}_{\theta}$ are Morita equivalent if and only if there exist a matrix $g \in SL_2 (\mathbb{Z}) $ such that $\theta' =g \theta$. 
\end{theorem}
In Section~\ref{modules} we will construct explicit bimodules realizing this equivalences. In what follows whenever we refer to a right $\mathcal{A}_{\theta}$-module (resp. left $\mathcal{A}_{\theta'}$-module, resp. 
$\mathcal{A}_{\theta'}$-$\mathcal{A}_{\theta}$-bimodule) we mean a  projective and finite-type
right $\mathcal{A}_{\theta}$-module (resp. left $\mathcal{A}_{\theta'}$-module, resp. 
$\mathcal{A}_{\theta'}$-$\mathcal{A}_{\theta}$-bimodule)

Given a irrational number $\theta$ a $\mathcal{A}_{\theta}$-$\mathcal{A}_{\theta}$-bimodule $E$ induces a Morita self equivalences of $\mathcal{A}_{\theta}$. We denote by $End_{Morita}(\mathcal{A}_\theta)$ the group of Morita self equivalence of $\mathcal{A}_\theta$. For example, given any positive integer  $n$ the free bimodule $\mathcal{A}_{\theta}^{n}$ induces a  a Morita self equivalences of $\mathcal{A}_{\theta}$. A Morita self equivalence defined via a free module is called a trivial  Morita self equivalence.

A Morita self equivalence of $\mathcal{A}_{\theta}$ given by a $\mathcal{A}_{\theta}$-$\mathcal{A}_{\theta}$-bimodule $E$ defines an endomorphism of $K_{0}(\mathcal{A}_\theta)$ via
\begin{eqnarray*}
\phi_E : [ P ] \mapsto [ P \otimes_{\mathcal{A}_\theta}  E ]
\end{eqnarray*} 
for $P$ a right projective finite rank  $\mathcal{A}_\theta$ module and $ [ P ] \in K_{0}(\mathcal{A}_\theta)$ is its $K$-theory class.

Via the map $\mathrm{rk}$ the endomorphism $\phi_M$ becomes multiplication by a real number. Thus we get a map:
\begin{eqnarray*}
\phi  : End_{Morita}(\mathcal{A}_\theta) \rightarrow 
\{ \alpha \in \mathbb{R} \, | \, \alpha \Gamma_\theta \subset \Gamma_\theta \}
\end{eqnarray*} 
this map turns out to be surjective. 

We can summarize the situation as follows (see \cite{Manin1}):
\begin{theorem}
Let $\theta\in\mathbb{R}$ be irrational. The following conditions are equivalent:
\begin{itemize} 
\item $\mathcal{A}_\theta $ has nontrivial Morita autoequivalences.
\item $\phi(End_{Morita}(\mathcal{A}_\theta)) \neq \mathbb{Z}$.
\item There exist  a matrix $g \in SL_2 (\mathbb{Z})$ such that 
\begin{eqnarray*}
\theta = g \theta .
\end{eqnarray*}
\item $\theta$ is a real quadratic irrationality: 
\begin{eqnarray*}
 [ \mathbb{Q}( \theta ) : \mathbb{Q}] = 2.
\end{eqnarray*}
\end{itemize}
\end{theorem}
If any of these equivalent conditions holds we say that the noncommutative torus $\mathbb{T}_\theta^2$ with algebra of smooth functions 
$\mathcal{A}_\theta$ is a {\it real multiplication noncommutative torus}. If $\mathbb{T}_\theta^2$ is a real multiplication noncommutative torus then
\begin{eqnarray*}
\phi(End_{Mor}(\mathcal{A}_\theta)) &=& \{ \alpha \in \mathbb{R} \, | \, \alpha \Gamma_\theta \subset \Gamma_\theta \} \\
&=& \mathbb{Z} + f \mathcal{O}_k 
\end{eqnarray*}
where $f\geq 1$ is an integer and $\mathcal{O}_k$ is the ring of integers of the real quadratic field $ k= \mathbb{Q}( \theta )$.

These results should be compared with the analogous results for elliptic curves leading to the theory of complex multiplication. The strong analogy suggests that  noncommutative tori may play a role in number theory similar to the role played by elliptic curves. In particular noncommutative tori with real multiplication could give the right geometric framework to attack the explicit class field theory problem for real quadratic fields (see \cite{Manin1}).  

\section{Heisenberg groups and their representations}
\label{SecHeis}

Various aspects of the theory of representations of Heisenberg groups arise \linebreak naturally when considering geometric constructions associated to noncommutative tori. This fact gives relations between noncommutative tori and elliptic curves through theta functions and plays a relevant role in the study of the arithmetic nature of related algebraic structures. In this section we sketch the parts of the theory of Heisenberg groups that are relevant in order to describe these results. We follow Mumford's Tata lectures \cite{tata3} which also recommend as a reference for the material in this section.

Let $G$ be a locally compact group lying in a central extension:
\begin{eqnarray*}
1  \,  \rightarrow  \,  \mathbb{C}^{*}_1 \,  \rightarrow  \, G   \,  {\rightarrow}  \, K  \,  \rightarrow  \, 0
\end{eqnarray*}
where $\mathbb{C}^{*}_1$ is the group of complex numbers of modulus $1$ and $K$ is a locally compact abelian group.  Assume that the exact sequence splits, so as a set  $G =\mathbb{C}^{*}_1 \times K$
and the group structure is given by:
\begin{eqnarray*}
(\lambda, x) (\mu, y ) &=& (\lambda \mu \psi(x, y)  , x + y ) 
\end{eqnarray*}
where $ \psi: K\times K \rightarrow  \mathbb{C}^{*}_1$ is a two-cocycle in $K$ with values in 
$\mathbb{C}^{*}_1$. The cocycle  $\psi$ induces a 
skew multiplicative pairing  
\begin{eqnarray*}
\label{pairing}
e: \; K\times K  & \longrightarrow & \; \mathbb{C}^{*}_1 \\ \nonumber 
(x,y)&\longmapsto& \frac{\psi(x, y)}{  \psi(y, x)}.
\end{eqnarray*}
This pairing defines a group morphism $\varphi: K \rightarrow  \widehat{K} $ from $K$ to its Pontrjagin dual given by 
$\phi(x)(y)=e(x,y)$.
\begin{definition} 
If $\varphi$ is a isomorphism we say that $G$ is a  \emph{Heisenberg group}.  
\end{definition}

For a Heisenberg group $G$ lying in a central extension as above we use the notation $G= \mathrm{Heis}(K)$. The main theorem about the representation of Heisenberg groups states that these kind of groups admit a unique normalized irreducible representations which can be realized in terms of a maximal isotropic subgroup of $K$. A subgroup $H$ of $K$ is called \emph{isotropic} if $e|_{H\times H} \equiv 1 $, this is equivalent to the existence of a section of $G$ over $H$:
\begin{eqnarray*}
\label{pairing2}
\sigma : \; K & \longrightarrow & G  \\ \nonumber 
x &\longmapsto& (\alpha(x),x).
\end{eqnarray*}
We say that a subgroup $H$ of $K$ is  \emph{maximal isotropic} if it is maximal with this property. A subgroup $H$ of $K$ is maximal isotropic if and only if $H = H^{\bot}$  where for $S \subset H$ we have 
\begin{eqnarray*}
\label{pairing2}
S^{\bot} &=& \{  x\in K \, | \, e(x,y)=1 \text{ for all } y \in S  \}.
\end{eqnarray*}

\begin{theorem}{\emph{(Stone, Von Neumann, Makey)}}
\label{SVNM}
Let $G$ be a Heisenberg group. Then 
\begin{itemize}
\item G has a unique irreducible unitary representation in which $\mathbb{C}^{*}_1$ acts by multiples of the identity. 
\item Given a maximal isotropic subgroup $H\subset K$ and a splitting $\sigma$ as above let 
$\mathcal{H }= \mathcal{H}_{H}$ be the space of measurable functions 
$f: K \rightarrow \mathbb{C}$ satisfying
\begin{enumerate}
\item $f(x+h)= \alpha(h)\psi(h,x)^{-1} f(x)$ for all $h \in H$.
\item $\int_{K / H} |f(x)|^2 dx  \, < \, \infty$.
\end{enumerate}
Then $G$ acts on $\mathcal{H}$ by
 \begin{eqnarray*}
\label{pairing2}
U_{(\lambda, y)} f (x) &=& \lambda \psi(x,y) f(x+y).
\end{eqnarray*}
and $\mathcal{H}$ is an irreducible unitary representation of $G$.
\end{itemize}
\end{theorem}
We call such representation a \emph{Heisenberg representation of $G$}.  The following theorem will be useful later:

\begin{theorem}{\emph{(Stone, Von Neumann, Makey)}}
\label{product}
Given two Heisenberg groups 
\begin{eqnarray*}
1  \,  \rightarrow  \,  \mathbb{C}^{*}_1 \,  \rightarrow  \, G_{i}   \,  {\rightarrow}  \, K_{i}  \,  \rightarrow  \, 0
\qquad i=1,2
\end{eqnarray*}
with Heisenberg representations $\mathcal{H}_{1}$ and $\mathcal{H}_{2}$ then 
\begin{eqnarray*}
1  \,  \rightarrow  \,  \mathbb{C}^{*}_1 \,  \rightarrow  
\, G_{1}\times G_{2} / \{(\lambda,\lambda^{-1}) | \lambda \in  \mathbb{C}^{*}_1  \}  
\,  {\rightarrow}  \, K_{1}\times K_{2}  \,  \rightarrow  \, 0
\end{eqnarray*}
is a Heisenberg group and its Heisenberg representation is $\mathcal{H}_{1}\hat{\otimes} \mathcal{H}_{2}$.
\end{theorem}

\subsection*{Real Heisenberg groups}
Let $K=\mathbb{R}^2$ and let $\varepsilon$ be a positive real number. 
We endow $G=K\times \mathbb{C}^{*}_1$ with the structure of a Heisenberg group defined by the  cocycle $\psi$ and the pairing $e$ given by 
\begin{eqnarray*}
 \psi(x, y) &=& \e \left(   \frac{1}{\varepsilon}\frac{(x_{1}y_{2}- y_{1}x_{2}) }{2} \right) \\
  e(x, y)  &=& \e \left(  \frac{1}{\varepsilon}(x_{1}y_{2}- y_{1}x_{2}) \right). 
\end{eqnarray*}
where $x=(x_{1},x_{2}), \, y =(y_{1},y_{2})\, \in K$.

If we choose as maximal isotropic subgroup $H=\{ x=(x_{1},x_{2})  \in K \, | \, x_{2} = 0 \}$ then the values of the functions in the corresponding Heisenberg representation  (Theorem~\ref{SVNM}) are determined by their values on  $\{ x=(x_{1},x_{2})  \in K \, | \, x_{1} = 0 \}$  and we may identify the space $\mathcal{H}_H$ with $L^{2}(\mathbb{R})$. The action of the $G$ is given by 
\begin{eqnarray*}
U_{(\lambda, y)} f (x) &=& 
\lambda \e \left(  \frac{1}{\varepsilon}  \left(  xy_{2}  + \frac{y_{1}y_{2} }{2} \right) \right)  f(x+y_{1}). 
\end{eqnarray*}
for $(\lambda, y) =(\lambda,(y_{1},y_{2}) )\in G$ and $f \in L^{2}(\mathbb{R})$. In particular we have 
\begin{eqnarray*}
U_{(1, (y_{1},0))} f (x) &=&   f(x+y_{1}). \\
U_{(1, (0,y_{2}))} f (x) &=&  \e \left(  \frac{1}{\varepsilon} xy_{2} \right)  f(x). 
\end{eqnarray*}
We will denote this Heisenber representation by $\mathcal{H}_{\varepsilon}$.
\medskip

For any $X\in Lie(G)$ and any Heisenberg representation $\mathcal{H}$ there is a dense subset of elements $f \in \mathcal{H }$ for which the limit 
\begin{eqnarray*}
\delta U _{X}(f) &=& \lim_{t\rightarrow 0}  \frac{U_{\exp(tX)}f   -f}{t} 
\end{eqnarray*}
exists. The above formula for $\delta U _{X}$ defines an unbounded operator on this set. An element $f \in \mathcal{H }$ is a smooth element for the representation $ \mathcal{H }$ of $G$ if 
\begin{eqnarray*}
\delta U _{X_{1}}\delta U _{X_{2}}\cdots \delta U _{X_{n}}(f)
\end{eqnarray*}
is well defined for any $n$ and any $X_{1},X_{2}, \dots X_{n} \in Lie(G)$. The set of smooth element of $\mathcal{H }$ is denoted by $\mathcal{H }_{\infty}$. We may realize $Lie(G)$ as an algebra of operators on $\mathcal{H }_{\infty}$. If choose a basis $\{A,B, C \}$ for the Lie algebra $Lie(G)$ such that
\begin{eqnarray*}
\exp(tA) =   (1,(t,0)), \quad \exp(tB) =   (1,(0,t)), \quad \exp(tC) =   (\e(t),(0,0)).
\end{eqnarray*} 
then a complex number $\tau\in \mathbb{C}$ with nonzero imaginary part  gives a decomposition of $Lie(G)\otimes \mathbb{C}$ into conjugate abelian complex subalgebras:
\begin{eqnarray*}
W_{\tau} &=& \langle \delta U _{A} - \tau \delta U _{B} \rangle \\
W_{\bar{\tau}} &=& \langle \delta U _{A} - \bar{\tau} \delta U _{B} \rangle .
\end{eqnarray*}

\begin{theorem}
\label{holom1}
Fix $\tau \in \mathbb{C}$ with $\Im(\tau) > 0$ then in any Heisenberg representation of $G$ there exists an element $f_{\tau}$, unique up to a scalar, such that $\delta U _{X}(f_{\tau})$ is defined and equal to $0$ for all $X\in W_{\tau}$. 
\end{theorem}

In the Heisenberg representation $\mathcal{H}_{\varepsilon}$ we have:
\begin{eqnarray*}
\delta U _{A}f(x) &=&   \frac{d}{dx}f(x) \\
\delta U _{B}f(x) &=&   \frac{2 \pi \imath x}{\varepsilon} f(x) \\
\delta U _{B}f(x) &=&   2 \pi \imath \, f(x) 
\end{eqnarray*} 
and $\mathcal{H }_{\varepsilon, \infty}$ is the Schwartz space $\mathcal{S}(\mathbb{R})$. The element $f_{\tau} $ in Theorem~\ref{holom1} is given by
\begin{eqnarray*}
\label{elemento}
f_{\tau} &=& \e \left(   \frac{1}{2\varepsilon} \tau x^{2}\right).
\end{eqnarray*}

\subsection{$\mathrm{Heis}(( \mathbb{Z}/c\mathbb{Z})^{2})$}
Let $c$ be a positive integer and let $ K=( \mathbb{Z}/c\mathbb{Z})^{2}$. We endow 
$G=K\times \mathbb{C}^{*}_1$ with the structure of a Heisenberg group defined by the  cocycle $\psi$ and the pairing $e$ given by 
 \begin{eqnarray*}
 \psi( (\left[  n_{1} \right], \left[n_{2} \right] ), (\left[  m_{1} \right], \left[m_{2} \right] ) )  
 &= &\e \left(  \frac{1}{2c}  \left( n_{1}m_{2}- m_{1}n_{2}  \right) \right) \\
 e( (\left[  n_{1} \right], \left[n_{2} \right] ), (\left[  m_{1} \right], \left[m_{2} \right] ) )  
 &= &\e \left(  \frac{1}{c}  \left( n_{1}m_{2}- m_{1}n_{2}  \right) \right)
\end{eqnarray*}
where $(\left[  n_{1} \right], \left[n_{2} \right] ), (\left[  m_{1} \right], \left[m_{2} \right] ) \, \in K$.

If we choose as maximal isotropic subgroup 
$H=\{ (\left[  n_{1} \right], \left[n_{2} \right])  \in K \, | \, \left[n_{2} \right] = 0 \}$ we may realize the Heisenberg representation as the action of $G$ on $C(\mathbb{Z}/c\mathbb{Z})$ given by 
\begin{eqnarray*}
U_{(\lambda, (\left[  m_{1} \right], \left[m_{2} \right] ))} \phi  (\left[ n \right] ) &=& 
\lambda \e \left(  \frac{1}{c}  \left(  nm_{2}  + \frac{m_{1}m_{2} }{2} \right) \right)  \phi( \left[ n +  m_{1}  \right] ). 
\end{eqnarray*}
for $(\lambda,(\left[  m_{1} \right], \left[m_{2} \right] )) \in G$ and $\phi\in C(\mathbb{Z}/c\mathbb{Z})$. In particular we have 
\begin{eqnarray*}
U_{(1, (\left[m_{1} \right] ,0))} \phi  (\left[ n \right] ) &=&    \phi  (\left[ n + m_{2} \right] ). \\
U_{(1, (0,\left[m_{2} \right] ))} \phi  (\left[ n \right] ) &=&  \e \left(  \frac{1}{c} nm_{2} \right)  f(\left[n \right] ). 
\end{eqnarray*}

\begin{remark}
This type of Heisenberg groups is related to \emph{algebraic Heisenberg groups}  or, more generally, \emph{Heisenberg group schemes}. These are given by central extensions of the form 
\begin{eqnarray*}
1  \,  \rightarrow  \,  \mathbb{G}_{m} \,  \rightarrow  \, \mathcal{G}   \,  {\rightarrow}  \,\mathcal{ K}  \,  \rightarrow  \, 0
\end{eqnarray*}
where $\mathcal{K}$ is a finite abelian group scheme aver a base field $k$. These groups arise in a natural way by considering ample line bundles on abelian varieties over the base field $k$. The corresponding Heisenberg representations can be realized as canonical actions of $\mathcal{G}$ on the spaces of sections of these bundles. The action of $Gal(\bar{k}/k)$ on the geometric points of $\mathcal{G}$ implies important algebraicity results about these representations. The abelian varieties that play a role 
in the constructions that follow are the elliptic curves  whose period lattice is spanned by the parameter 
$\tau$ which defines the complex structure on the noncommutative torus.  
\end{remark}

\section{Heisenberg modules over noncommutative tori with real multiplication}
\label{modules} 
Let $\theta \in \mathbb{R}$ be a quadratic irrationality and let 
\begin{eqnarray*}
g&=& \left(
\begin{array}{cc}
  a & b \\
  c & d \\
\end{array}
\right)\in SL_2 (\mathbb{Z})
\end{eqnarray*}
be a matrix fixing $\theta$. In this section we describe the construction of  a $\mathcal{A}_{\theta}$-$\mathcal{A}_\theta$-bimodule $E_g$ whose isomorphism class gives a Morita self equivalence of $\mathcal{A}_{ \theta}$. In what follows we assume that $c$ and $c\theta +d$ are positive. 
Let
\begin{eqnarray*}
\varepsilon &=&  \frac{c\theta +d}{c}
\end{eqnarray*}
and consider the following operators on the Schwartz space $\mathcal{S}(\mathbb{R})$:
\begin{eqnarray*}
(\check{U}f)(x) &=& f(x- \varepsilon)  \\
(\check{V}f)(x) &=& \e(x) f(x) \\
(\hat{U}f)(x) &=& f \left( x - \frac{ 1}{c} \right)  \\
(\hat{V}f)(x) &=& \e \left(\frac{x}{c \varepsilon} \right)  f(x).
\end{eqnarray*}
Note that each pair of operators corresponds to the Heisenberg group action of two generators of 
$\mathbb{R}^{2}$ where as above we identify the Schwartz space $\mathcal{S}(\mathbb{R})$ with the set of smooth elements of the Heisenberg representation $\mathcal{H}_{\varepsilon}$ of $\mathrm{Heis}(\mathbb{R}^{2})$.

We consider also the following operators on $C(\mathbb{Z}/c\mathbb{Z})$
\begin{eqnarray*}
(\check{u} \phi)(\left[ n \right]) &=& \phi(\left[ n - 1 \right])  \\
(\check{v} \phi) (\left[ n \right])  &=& \bar{\e} \left(\frac{ d n}{c} \right)  \phi (\left[ n \right]) \\
(\hat{u} \phi)(\left[ n \right]) &=& \phi(\left[ n - a \right])  \\
(\hat{v} \phi) (\left[ n \right])  &=& \bar{\e} \left(\frac{ n}{c} \right)  \phi (\left[ n \right]) .
\end{eqnarray*}
Since both $a$ and $d$ are prime relative to $c$ each pair of operators corresponds to the Heisenberg group action of two generators of 
$(\mathbb{Z}/c\mathbb{Z})^{2}$ on the Heisenberg representation $C(\mathbb{Z}/c\mathbb{Z})$ of  
$\mathrm{Heis}((\mathbb{Z}/c\mathbb{Z})^{2})$.

Taking into account the commutation relations satisfied between each of the above pairs of operators and the fact that $g \theta = \theta$ we see that  the space 
\begin{eqnarray*}
E_{g} &=& \mathcal{S}(\mathbb{R})\otimes  C(\mathbb{Z}/c\mathbb{Z})
\end{eqnarray*}
becomes a $\mathcal{A}_{\theta}$-$\mathcal{A}_\theta$-bimodule by defining:
\begin{eqnarray*}
(f\otimes \phi )U  &=& (\check{U} \otimes  \check{u} ) (f\otimes \phi)   \\
(f\otimes \phi )V &=& (\check{V} \otimes  \check{v} ) (f\otimes \phi)   \\
U(f\otimes \phi )  &=& (\hat{U} \otimes  \hat{u} ) (f\otimes \phi)   \\
V(f\otimes \phi )  &=& (\hat{V} \otimes  \hat{v} ) (f\otimes \phi)   
\end{eqnarray*}
where $f\in  \mathcal{S}(\mathbb{R})$ and $\phi \in C(\mathbb{Z}/c\mathbb{Z})$.

\begin{theorem}\emph{(Connes \cite{Connes1})}
With the above bimodule structure $E_{g} $ is finite-type and projective both as a right  $\mathcal{A}_\theta$-module  and as a left  $\mathcal{A}_\theta$-module. Considered as a right module we have 
$\mathrm{rk}(E_{g})= c \theta +d $. 
The left action of $\mathcal{A}_\theta$ gives an identification: 
\begin{eqnarray*}
\mathrm{End}_{\mathcal{A}_\theta}(E_g) \simeq \mathcal{A}_{\theta}.
\end{eqnarray*}
\end{theorem}
We refer to this kind of modules as \emph{Heisenberg modules}.

\medskip
Taking into account the $\mathcal{A}_\theta$-$\mathcal{A}_\theta$-bimodule structure of $E_{g}$ we may consider the tensor product $E_{g}\otimes_{\mathcal{A}_\theta}E_{g}$. This is one of the main consequences of the real multiplication condition. There is a natural identification (see \cite{DiS, PoS}):
\begin{eqnarray*}
\label{producto1}
E_{g}\otimes_{\mathcal{A}_\theta}E_{g}   &\simeq & 
E_{g^{2}} .
\end{eqnarray*}

To see this consider first the completed tensor product over $\mathbb{C}$ of the space $E_{g}$ with it self:
\begin{eqnarray*}
E_{g} \hat{\otimes} E_{g}  &=& 
\left[  \mathcal{S}(\mathbb{R}) \otimes  C(\mathbb{Z}/c\mathbb{Z})  \right]  \hat{\otimes} 
\left[\mathcal{S}(\mathbb{R}) \otimes  C(\mathbb{Z}/c\mathbb{Z})   \right]       \\
&=& \left[ \mathcal{S}(\mathbb{R})  \hat{\otimes}  \mathcal{S}(\mathbb{R}) \right]  \otimes  
\left[ C(\mathbb{Z}/c\mathbb{Z})  \otimes  C(\mathbb{Z}/c\mathbb{Z}) \right]  \\
&=& \left[ \mathcal{S}(\mathbb{R} \times \mathbb{R})  \right]  \otimes  
\left[ C(\mathbb{Z}/c\mathbb{Z} \times \mathbb{Z}/c\mathbb{Z}) \right]  .
\end{eqnarray*}

The space $\mathcal{S}(\mathbb{R} \times \mathbb{R})$ is the space of smooth elements of the Heisenberg representation of $\mathrm{Heis}(\mathbb{R}^{4})$ obtained as a product of the Heisenberg representations of $\mathrm{Heis}(\mathbb{R}^{2})$. Likewise $C(\mathbb{Z}/c\mathbb{Z} \times \mathbb{Z}/c\mathbb{Z})$ is the Heisenberg representation of $\mathrm{Heis}((\mathbb{Z}/c\mathbb{Z})^{4})$ obtained as a product of the Heisenberg representations of $\mathrm{Heis}((\mathbb{Z}/c\mathbb{Z})^{2})$. 

To pass from $E_{g} \hat{\otimes} E_{g}$ to $E_{g}  \otimes_{\mathcal{A}_{\theta}} E_{g}$ we have to quotient
$E_{g} \hat{\otimes} E_{g}$ by the space spanned by the relations:
\begin{eqnarray*}
\left[ (f\otimes \phi )U \right]  \hat{\otimes} \left[ g \otimes \omega \right]  
&=& \left[ f\otimes \phi  \right]  \hat{\otimes} \left[ U (g \otimes \omega ) \right]    \\
\left[ (f\otimes \phi )V \right]  \hat{\otimes} \left[ g \otimes \omega \right]  
&=& \left[ f\otimes \phi  \right]  \hat{\otimes} \left[ V (g \otimes \omega ) \right]  
\end{eqnarray*}
where $f,g\in  \mathcal{S}(\mathbb{R})$ and $\phi, \omega \in C(\mathbb{Z}/c\mathbb{Z})$.
 
 At the level of the Heisenberg representations involved this amounts to restrict to the subspaces of 
 $\mathcal{S}(\mathbb{R} \times \mathbb{R})$ and $C(\mathbb{Z}/c\mathbb{Z} \times \mathbb{Z}/c\mathbb{Z})$ which are invariant under the action of the subgroups of  $\mathrm{Heis}(\mathbb{R}^{4})$ and  $\mathrm{Heis}((\mathbb{Z}/c\mathbb{Z})^{4})$  generated by the elements giving these relations.
 
The corresponding space of invariant elements  in $\mathcal{S}(\mathbb{R} \times \mathbb{R})$ is canonically isomorphic to the space  $\mathcal{S}(\mathbb{R})$ of smooth elements of the Heisenberg representation of $\mathrm{Heis}(\mathbb{R}^{2})$ with $c\varepsilon^{2}/(a+d)$ playng the role of $\varepsilon$. In $C((\mathbb{Z}/c\mathbb{Z})\times(\mathbb{Z}/c\mathbb{Z}))$ the corresponding invariant subspace is canonically isomorphic to the Heisenber representation $C(\mathbb{Z}/c(a+d)\mathbb{Z})$ of $\mathrm{Heis}((\mathbb{Z}/c(a+d)\mathbb{Z})^{2})$. Thus we get:
\begin{eqnarray*}
E_{g} \otimes_{\mathcal{A}_{\theta}} E_{g} &\simeq& \mathcal{S}(\mathbb{R}) \otimes  C(\mathbb{Z}/c(a+d)\mathbb{Z}) \\
&=& E_{g^{2}}  
\end{eqnarray*}
The compatibility of the module structures in this isomorphism is implied by the the compatibility of the Heisenberg representations involved.

In a similar manner one may obtain isomorphisms:  
\begin{eqnarray*}
\underbrace{E_{g}\otimes_{\mathcal{A}_ {\theta}} \dots
\otimes_{\mathcal{A}_ {\theta}}  E_{g}}_{n} 
\simeq  E_{g^{n}}.
\end{eqnarray*}

\section{Some rings associated to noncommutative tori with real multiplication}

Noncommutative tori may be considered as noncommutative projective varieties. In noncommutative algebraic geometry varieties are defined in terms of categories which play the role of appropriate categories of sheaves on them (see \cite{StaffordVdBergh}).  In \cite{Po} Polishchuk  analyzed real multiplication  noncommutative tori from this point of view. Given a real quadratic irrationality $\theta$ and a complex structure $\delta_{\tau}$ on $\mathcal{A}_{\theta}$  Polishchuk constructs a homogeneous coordinate ring associated $\mathbb{T}_{\theta}^2$ and $\delta_{\tau}$.  

During the rest of the section $g$ will denote a matrix in $ SL_2 (\mathbb{Z})$ fixing a quadratic irrationality $\theta$. We will denote the elements of the powers of this matrix by 
\begin{eqnarray*}
g^{n}= \left(
\begin{array}{cc}
  a_{n} & b_{n} \\
  c_{n} & d_{n} \\
\end{array}
\right),  \qquad n > 0. 
\end{eqnarray*}
Given a Heisenberg $\mathcal{A}_{\theta}$-$\mathcal{A}_{\theta}$-bimodule $E_{g}$ we use the complex structure  
$\delta_{\tau}$ on $\mathcal{A}_{\theta}$ to single out a finite dimensional subspace in each of the graded pieces of
\begin{eqnarray*}
\mathcal{E}_{g} &=&\bigoplus_{n\geq 0} \underbrace{E_{g}\otimes_{\mathcal{A}_ {\theta}} \dots
\otimes_{\mathcal{A}_ {\theta}}  E_{g}}_{n} \\
 &=& \bigoplus_{n\geq 0}E_{g^{n}}  .
\end{eqnarray*}
This should be done in a way compatible with the product structure of $\mathcal{E}_{g}$.

Given a Heisenberg $\mathcal{A}_{\theta}$-$\mathcal{A}_{\theta}$-bimodule $E_{g}$ we may define a connection on $E_{g}$ by: 
\begin{eqnarray*}
(\nabla_{1} f \otimes \phi)(x,\left[ n \right]) &=& 2 \pi \imath \left(   \frac{x}{\varepsilon} \right) (f \otimes \phi)(x,  \left[ n \right])  \\
(\nabla_{2}f \otimes \phi)(x,\left[ n \right]) &=&\frac{ d}{dx}(f \otimes \phi)(x,  \left[ n \right]) .
\end{eqnarray*}
Connections of this kind were studied in \cite{ConnesRieffel} in the context of Yang Mils theory for noncommutative tori. Note that this connection corresponds to the action of the Lie algebra of $\mathrm{Heis}(\mathbb{R}^{2})$ on the left factor of $\mathcal{S}(\mathbb{R}) \otimes  C(\mathbb{Z}/c\mathbb{Z})$ given by the derivations 
$\delta U_{A}$ and $\delta U_{B}$. Once we choose a complex parameter $\tau \in \mathbb{C} \setminus \mathbb{R}$ giving a complex structure on $\mathcal{A}_{\theta}$ the corresponding decomposition of the complexified Lie algebra singles out the element $f_{\tau}$ 
(Theorem~\ref{holom1}) thus it is natural to consider the space:
\begin{eqnarray*}
R_{g} &=& \{ f_{\tau}\otimes \phi \in E_{g} \, | \, \phi \in C(\mathbb{Z}/c\mathbb{Z}) \} \\
 &=& \Big\{ \e \left(   \frac{1}{2\varepsilon} \tau x^{2}\right) \otimes \phi \in E_{g} \, | \, \phi \in C(\mathbb{Z}/c\mathbb{Z})   \Big\}.
\end{eqnarray*} 
This are the  the spaces of  the holomorphic vectors considered in \cite{Schwarz1,DiS, PoS,Po}.

We denote by $f_{\tau,n}$ the corresponding element on the left factor of $E_{g^{n}}$ and let 
\begin{eqnarray*}
R_{g^{n}} &=& \Big\{ f_{\tau,n}  \otimes \phi \in E_{g^{n}} \, | \, \phi \in C(\mathbb{Z}/c_{n}\mathbb{Z})   \Big\}.
\end{eqnarray*} 
Following \cite{Po} we define the  {\it homogeneous coordinate ring for the noncommutative torus 
$\mathbb{T}_{\theta}^2$  with complex structure $\delta_{\tau}$} by:
\begin{eqnarray*}
B_g(\theta,\tau) &=& \bigoplus_{n\geq 0}R_{g^{n}} .
\end{eqnarray*}
The following result characterizes some structural properties of $B_g(\theta,\tau)$ in terms of the matrix elements of $g$:

\begin{theorem}\emph{(\cite{Po} Theorem 3.5)}
\label{poli2}
Assume $g\in SL_2 (\mathbb{Z})$ has positive real eigenvalues: 
\begin{enumerate}
\item If $c \geq a+ d$ then $B_g(\theta,\tau)$ is generated over $\mathbb{C}$ by $R_{g}$.
\item If $c \geq a+ d+1$ then $B_g(\theta,\tau)$ is a quadratic algebra.
\item \label{la3} If $c \geq a+ d+2$ then $B_g(\theta,\tau)$ is a Koszul algebra.
\end{enumerate}
\end{theorem}
\medskip

Let $X_{\tau}$ be the elliptic curve with complex points $\mathbb{C}/ (\mathbb{Z} \oplus\tau \mathbb{Z})$. Taking into account the remarks at the end of section Section~\ref{SecHeis} it is possible to realize each space $R_{g^{n}}$ as the space of sections of a line bundle over $X_{\tau}$. For this we consider the matrix coefficients obtained by pairing $f_{\tau,n}$ with functionals in the distribution completion of the Heisenberg representation  which are invariant under the action of elements in $\mathrm{Heis}(\mathbb{R}^{2})$ corresponding to a lattice in $\mathbb{R}^{2}$ associated to $g^{n}$. This matrix coefficients correspond to theta functions with rational characteristics which form a basis for the space of sections of the corresponding line bundle over $X_{\tau}$. In these bases the structure constants for the product of $B_{g}(\theta,\tau)$ have the form
\begin{eqnarray*}
\vartheta_{r}(l \tau)
\end{eqnarray*}
where $\l\in\mathbb{Z}$ and $\vartheta_{r}(l \tau)$ is the theta constant with rational characteristic 
$r\in \mathbb{Q}$
defined by the series 
\begin{eqnarray*}
\label{thetaseries}
\vartheta_{r}(l \tau) = \sum_{n \in \mathbb{Z}} \exp [\pi \imath (n+r)^{2} l \tau].
\end{eqnarray*}

This fact has the following consequence:
\begin{theorem}{ \emph{(\cite{yo})}}
Let $\theta \in \mathbb{R}$ be a quadratic irrationality fixed by a matrix $g\in SL_2 (\mathbb{Z})$ and assume $c \geq a+ d+2$. Let $k$ be the minimal field of definition of the elliptic curve $X_{\tau}$.
Then the algebra $B_{g}(\theta,\tau)$ admits a rational presentation over a finite algebraic extension of $k$.
\end{theorem}

\begin{remark}
Analogous results hold for the rings of quantum theta functions considered in \cite{Masha}.
These rings correspond to Segre squares of  the homogeneous coordinate rings $B_{g}(\theta,\tau)$ and can be analyzed in terms of the Heisenberg modules involved in their construction.
\end{remark}

\bibliographystyle{amsalpha}

\begin{thebibliography}{A}

\bibitem[1]{Connes1} A. Connes, \textit{$C^{*}$-alg\`{e}bres et g\'{e}om\'{e}trie diff\'{e}rentielle.} C. R. Ac. Sci. Paris, t. \textbf{290}   (1980) 599-604.

\bibitem[2]{Connes3} A. Connes, \textit{Noncommutative geometry.} Academic Press, 1994.

\bibitem[3]{ConnesRieffel} A. Connes, M. Rieffel, \textit{Yang-Mills for noncommutative two-tori. Operator algebras and mathematical physics} Contemp. Math., \textbf{62}, Amer. Math. Soc. Providence, RI (1987) 237-266. 
   
\bibitem[4]{Davidson} K. R. Davidson, \textit{$C^{*}$-algebras by example.}
Fields Institute Monographs, 6. American Mathematical Society, Providence, RI, (1996).

\bibitem[5]{DiS} M. Dieng, A. Schwarz, \textit{Differential and complex geometry of two-dimensional noncommutative tori.} Lett. Math. Phys.  \textbf{61}  (2002) 263-270.

\bibitem[6]{Effros} E. G. Effros, F. Hahn, \textit{Locally compact transformation groups and 
$C^*$-algebras.}  Bull. Amer. Math. Soc.  \textbf{73}  (1967) 222-226. 

\bibitem[7]{Manin1} Y. Manin \textit{Real multiplication and noncommutative geometry}. In: ``The legacy of Niels Henrik Abel", Springer Verlag, Berlin (2004).

\bibitem[8]{tata3} D. Mumford \textit{Tata Lectures on Theta III.} With the collaboration of M. Nori and P. Norman. Progr. Math. 97. Birkh\"auser, Boston (1991).

\bibitem[9]{yo} J. Plazas \textit{Arithmetic structures on noncommutative tori with real multiplication}   
Preprint: ArXiv math.QA/0610127
    
\bibitem[10]{Po} A. Polishchuk, \textit{Noncommutative two-tori with real multiplication as noncommutative projective varieties.} Journal of Geometry and Physics \textbf{50} (2004) 162-187.

\bibitem[11]{Polishchuk2} A. Polishchuk \textit{Classification of holomorphic vector bundles on noncommutative two-tori.}  Doc. Math.  \textbf{9}  (2004) 163-181.

\bibitem[12]{PoS} A. Polishchuk, A. Schwarz. \textit{Categories of holomorphic bundles on noncommutative two-tori}. Comm. Math. Phys. \textbf{236} (2003) 135-159.

\bibitem[13]{rieffel1} M. A. Rieffel, \textit{$C^{*}$-algebras associated with irrational rotations.} Pacific J. Math.  \textbf{93} (1981) 415-429.

\bibitem[14]{rieffel2} Rieffel, M. A. \textit{Noncommutative tori - a case 
study of noncommutative differentiable manifolds}. Geometric and topological
invariants of elliptic operators, Contemp. Math., 105, Amer. Math. Soc., Providence, RI, (1990) 191-211.

\bibitem[15]{riff3} Rieffel, M.A. \emph{The cancellation theorem for 
projective modules over irrational rotation $C^{*}$-algebras}. Proc. London
Math. Soc. (3) 47, 285-302 (1983).

\bibitem[16]{Schwarz1} A. Schwarz, \textit{Theta functions on noncommutative tori.} Lett. Math. Phys.  \textbf{58}  (2001) 81-90.

\bibitem[17]{StaffordVdBergh} J. T. Stafford, M. van den Bergh, \textit{Noncommutative curves and noncommutative surfaces.}  Bull. Amer. Math. Soc. (N.S.)  \textbf{38}  (2001) 171-216.

\bibitem[18]{Masha} M. Vlasenko, \textit{The graded ring of quantum theta functions for noncommutative torus with real multiplication}.  Int. Math. Res. Not. Art. ID 15825 (2006).


\end{thebibliography}

\end{document}